\documentclass[twoside,11pt,leqno]{article}

\usepackage{amssymb}
\usepackage{amsmath}
\usepackage{amsfonts}
\usepackage{enumerate}

\usepackage{amsthm}
\newtheorem{defn}{Definition}[section]
\newtheorem{teo}{Theorem}[section]

\newtheorem{lem}{Lemma}[section]
\newtheorem{pro}{Proposition}[section]
\newtheorem{cor}{Corollary}[section]
\newtheorem{rem}{Remark}[section]

\let\oldproofname=\proofname
\renewcommand{\proofname}{\rm\bf{\oldproofname}}

\usepackage[all]{xy}

\begin{document}

\title{Eventually constant intertwining linear maps between complete locally convex spaces}
\maketitle

 {\bf Carlos Bosch}  \\
 Departamento de Matem\'{a}ticas,  ITAM,  R\'{\i}o Hondo \#1, Col. Progreso Tizap\'{a}n, Mexico, DF, 01080,  Mexico  \\ 
bosch@itam.mx \\

{\bf C\'{e}sar L.  Garc\'{i}a}  \\
     Departamento de Matem\'{a}ticas,  ITAM,  R\'{\i}o Hondo \#1, Col. Progreso Tizap\'{a}n, Mexico, DF, 01080,   Mexico\\
 clgarcia@itam.mx   \\
 
 {\bf Thomas E. Gilsdorf}  \\
      Department of Mathematics,  Central Michigan University,  Mt. Pleasant, MI  48859    USA \\
     gilsd1te@cmich.edu    \\
     
 {\bf Claudia G\'{o}mez - Wulschner}   \\
Departamento de Matem\'{a}ticas,  ITAM,  R\'{\i}o Hondo \#1, Col. Progreso Tizap\'{a}n, Mexico, DF, 01080, Mexico   \\
claudiag@itam.mx   \\

  {\bf Rigoberto Vera}    \\
     Departamento de Matem\'{a}ticas,  ITAM,  R\'{\i}o Hondo \#1, Col. Progreso Tizap\'{a}n, Mexico, DF, 01080, Mexico   \\
     rveram@itam.mx.   \


\begin{abstract}    Starting from Sinclair's 1976 work [6] on automatic continuity of linear operators on Banach spaces, we prove that sequences of intertwining continuous linear maps are eventually constant with respect to the separating space of a fixed linear map.  Our proof uses a gliding hump argument.    We also consider aspects of continuity of linear functions  between locally convex spaces and prove that  such  a linear function $T$ from  the locally convex space $X$ to the locally convex space $Y$ is continuous whenever the separating space $G(T)$ is the zero vector in $Y$ and for which $X$ and $Y$ satisfy conditions for a closed graph theorem.    \\

\end{abstract}


\section{Introduction}

\noindent  From  A. Sinclair's fundamental work \cite{SC}, one of the main result with respect to linear maps  between Banach spaces is Lemma 1.6 (\cite{SC}; p. 11), which gives conditions under which a sequence of intertwining linear maps is eventually constant: \\

\noindent {\bf Lemma} (\cite{SC}).  {\it Suppose $X$ and $Y$ are Banach spaces and let $(T_{n})$ and $(R_{n})$ be sequences of continuous linear operators on $X$ and $Y$, respectively.  If $S$ is a linear operator from $X$ to $Y$ that satisfies $S T_{n} = R_{n}S$ for all $n\in \mathbb{N}$, then there is an integer $N\in \mathbb{N}$ such that} $\forall n\geq N$

$$   \overline{R_{1} \cdots R_{n} \mathcal{G}} =  \overline{R_{1} \cdots R_{N} \mathcal{G}},   $$

\noindent {\it where $\mathcal{G}$ is the separating space of } $S$:
$$ \mathcal{G} =  \{ y\in Y : \mbox{there is a sequence } (x_{n})\rightarrow 0 \mbox{ in } X \ni (Sx_{n}) \rightarrow y \mbox{ in } Y\}.  $$

\noindent The proof utilizes a gliding hump argument.  See \cite{CS} for example.  The importance of the result on eventually constant sequences of intertwining maps is in its application to results on automatic continuity.  Some instances include  results on automatic continuity in  \cite{SC} for Banach spaces (see also the references therein), in  \cite{NA}, for a class of locally convex spaces, and within the context of generalized local linear operators, and in  \cite{AT}, in the context of convex bornological spaces where the results are given with regard to bounded linear maps.   In this paper, we define the separating space for general locally convex spaces and prove  that Sinclair's  result  holds in for sequences of intertwining linear maps between complete locally convex spaces.  This opens the door to applications in general locally convex spaces for which only completeness is needed.    Like in the proof of the Banach space result, our proof utilizes a gliding hump argument. Similar generalizations to complete locally convex spaces or related spaces have been obtained recently, such as in  \cite{HQ}, and  \cite{QH}.        

\noindent  The notation here is as follows.  Let $(X,\tau)$  be a Hausdorff locally convex topological space (lcs)  over $\mathbb{K} = \mathbb{R} \mbox{ or } \mathbb{C}$.  Denote by  $\{\rho_{j}\}_{J}$  a family of seminorms that generate $\tau$ on $X$.   Let $Z_{\rho_{j}} =\{x\in X : \rho_{j}(x)=0\}$. We denote the zero vector of $X$ by $\hat{0}$.   Observe that $Z_{\rho_{j}}$ is a closed linear subspace of $X$ and $\bigcap_{J}Z_{\rho_{j}}=\{\hat{0}\}$.  We denote a linear subspace $M$ of a linear space $L$ by $M < L$.  \\

\noindent  Let $\,(Y,\sigma)$ be another lcs, with $H = \{\sigma_{i}\}_{I}$ a family of seminorms generating $\sigma$.  Denote by $Z_{\sigma_{i}}$ the corresponding $\sigma$- closed linear subspace of $Y$ at which   $\sigma_{i} = 0$ in $Y$. We denote  the zero of $Y$ by $\tilde{0}$.    \\

\noindent Let $\,T:(X,\tau) \rightarrow (Y,\sigma)\,$ be a linear map, not necessarily continuous,  and notice that $T(Z_{\rho_{j}})< Z_{\sigma_{i}}\,$ need not hold.  

\noindent  As necessary, we may also include a third lcs $(W,\lambda)$, with corresponding family of seminorms $\{\lambda_{k}\}_{K}\,$ that generates $\lambda$.  \\

\noindent  All other unspecified notation is standard, e.g., as in  \cite{NB}.    \\
 

\noindent  The definition below, inspired by     \cite{SC},  page 8,  is defined  for locally convex spaces, using seminorms.  

\begin{defn}
For each $\,j\in J\,$\, and each $\,i\in I\,$, let 

$$G_{ji}(T)=\{y\in Y :  \exists \{x_{n}\}_{\mathbb{N}}\stackrel{\rho_{j}}{\rightarrow} \hat{0}\in X \mbox{ and } \{Tx_{n}\}_{\mathbb{N}}\stackrel{\sigma_{i}}{\rightarrow} y\}.  $$

\noindent We then define 
$$\,G(T)=\bigcap_{J\times I}G_{ji}(T)=\{y\in Y :   \exists \{x_{n}\}_{\mathbb{N}}\stackrel{\tau}{\rightarrow} \hat{0}\in X  \mbox{ and } \{Tx_{n}\}_{\mathbb{N}}\stackrel{\sigma}{\rightarrow} y\}. $$ 
\end{defn}

\noindent We have $\,G_{ji}(\alpha T)=\alpha G_{ji}(T)\,$ \hspace{2mm} $\,\forall\, \alpha \in \mathbb{K}\,$.

\section{Preliminary results.}  

\begin{pro}
 For each $\,(j,i)\in J\times I\,$ $\,G_{ji}(T)\,$ satisfies:
\begin{enumerate}

\item  $\,G_{ji}(T)\,$ is a linear subspace of $Y$.
\item $\,G_{ji}(T)\,$ is $\,\sigma_{i}$-closed.
\item $\,G_{ji}(T)\,$ is $\,\sigma$-closed.
\item $\,Z_{\sigma_{i}}<G_{ji}(T)\,$.   
\item  $\,T(Z_{\rho_{j}})<G_{ji}(T)\,$. 
\end{enumerate}
\end{pro}


\noindent {\it Proof} 
 1. Let $\,y,  y'\in G_{ji}(T)\,$\ and $\,\lambda\in \mathbb{K}$. There exist sequences  $\,\{x_{n}\}_{\mathbb{N}}\stackrel{\rho_{j}}{\rightarrow} \hat{0}\,$\ and  $\,\{x_{n}'\}_{\mathbb{N}}\stackrel{\rho_{j}}{\rightarrow} \hat{0}\,$\, in $X$  such that  $\,\{T(x_{n})\}_{\Bbb{N}}\stackrel{\sigma_{i}}{\rightarrow} y\,$ and  $\,\{T(x_{n}')\}_{\Bbb{N}}\stackrel{\sigma_{i}}{\rightarrow} y'\,$.  From this $\,\{x_{n}+\lambda x_{n}'\}_{\mathbb{N}}\stackrel{\rho_{j}}{\rightarrow} \hat{0}\,$ and  $\,\{T(x_{n}+\lambda x_{n}')\}_{\mathbb{N}}\stackrel{\sigma_{i}}{\rightarrow} y+\lambda y'\,$.   Hence,  $\,y+\lambda y'\in G_{ji}(T)\,$.    \\


2.  Let  $\,\{y_{k}\}_{\mathbb{N}}\subset G_{ji}(T)\,$ such that  $\,\{y_{k}\}_{\mathbb{N}}\stackrel{\sigma_{i}}{\rightarrow} y\in Y\,$.   For each  $\,y_{k}\,$\, there exists \, $\,\{x_{kn}\}_{n\in \mathbb{N}}\stackrel{\rho_{j}}{\rightarrow} \hat{0}\,$ and $\,\{T(x_{kn})\}_{n\in\mathbb{N}}\stackrel{\sigma_{i}}{\rightarrow} y_{k}\,$ as $\,n\rightarrow \infty\,$.  For each  $\,\epsilon_{k}=\frac{1}{k}\,$ let $\,N_{k}\in \mathbb{N}\,$ such that $\,\rho_{j}(x_{kn})<\epsilon_{k}\,$ and  $\,\sigma_{i}(T(x_{kn})-y)<\epsilon_{k}\,$,  $\,\forall\, n\geq N_{k}\,$.  Let  $\,\{x_{k}=x_{kN_{k}}\}\subset X\,$.  \\

\noindent {\it Claim 1}:\,  $\,\{x_{k}\}_{\mathbb{N}}\stackrel{\rho_{j}}{\rightarrow} \hat{0}\,$.  \
 
 \noindent To see this, let $\,\epsilon >0\,$ and consider  $\,K\in \mathbb{N}\,$ such that  $\,\frac{1}{K}<\epsilon\,$:\, $\,\rho_{j}(x_{kN_{k}})<\frac{1}{K}<\epsilon\,$ if  $\,k\geq K\,$.   \\

\noindent {\it Claim 2}:\,  $\,\{T(x_{k})\}_{\mathbb{N}}\stackrel{\sigma_{i}}{\rightarrow} y\,$.    \

\noindent  For this, let  $\,\epsilon >0\,$ and $\,K\in \mathbb{N}\,$ such that \, $\,\frac{1}{K}< \frac{\epsilon}{2}\,$.  Then $\,\forall\, k\geq K\,$:  

$$  \sigma_{i}(T(x_{k})-y)\leq \sigma_{i}(T(x_{k})-y_{k})+\sigma_{i}(y_{k}-y)<\frac{1}{K}+\frac{\epsilon}{2}<\epsilon.   $$

\noindent    Thus, $\,y\in G_{ji}(T)\,$, which proves  $\,G_{ji}\,$\, is sequentially  $\,\sigma_{i}$-closed, and we also conclude $\,\sigma_{i}$-closed.   \\

\noindent 3.  Consider the net $\,\{y_{\lambda}\}_{\Lambda}\subset G_{ji}(T)\,$ such that $\,\{y_{\lambda}\}_{\Lambda}\stackrel{\sigma}{\rightarrow} y\in Y\,$.  Because the topology generated by $\,\sigma_{i}\,$ is contained in the topology $\sigma\,$, we have $\,\{y_{\lambda}\}_{\Lambda}\stackrel{\sigma_{i}}{\rightarrow} y\in Y\,$.  By item 2,  $\,y\in G_{ji}(T)\,$.    \\

\noindent 4. Let  $\,y\in Z_{\sigma_{i}}\,$.   Then  $\,\sigma_{i}(y)=0\,$.  Take the sequence $\,\{\hat{0}, \hat{0}, \hat{0}, ... \}\,$, in $X$,  which  converges to $\,\hat{0}\,$.  We show that $\,\{T(\hat{0}), T(\hat{0}), T(\hat{0}), ... \}=\{\tilde{0}, \tilde{0}, \tilde{0}, ... \} \stackrel{\sigma_{i}}{\rightarrow} y\,$. Indeed, $\sigma_{i}(T(\hat{0})-y)\leq \sigma_{i}(\tilde{0})+\sigma_{i}(y)=0\,$.  It follows that $\,y\in G_{ji}(T)\,$.    \\

\noindent 5.  Let   $\,x\in Z_{\rho_{j}}\,$.  Then we have $\,\rho_{j}(x)=0\,$.  This time, we take the sequence $\{x, x, x, ... \}$ in $X$, which is $\rho_{j}$- convergent  to $\hat{0}$.  Put  $\,y=T(x)\,$.  Then $\,\{T(x), T(x), T(x), ... \}\stackrel{\sigma_{i}}{\rightarrow} y\,$.  We conclude that  $\,T(x)=y\in G_{ji}(T)\,$.   \/ \/ $\Box$  \\

\noindent We observe that, given a linear $T: X\rightarrow Y$ that need not be continuous, the subspaces $\,T(Z_{\rho_{j}})\,$ and   $\,Z_{\sigma_{i}}\,$ have a nontrivial intersection.  This occurs despite the fact that $\,T(Z_{\rho_{j}})\,$ might not be a subspace of  $\,Z_{\sigma_{i}}\,$.   On the other hand,  not every linear map will send the linear subspace $Z_{\rho_{j}}\,$ of the domain into the subspace $\,Z_{\sigma_{i}}\,$ of the codomain.  However, item 5 in the above proposition is somewhat surprising because it indicates that all linear maps will send $\,Z_{\rho_{j}}\,$ into $\,G_{ji}(T)\,$, even though the general size of $\,Z_{\sigma_{i}}< G_{ji}(T)\,$ is not much larger.  We will see shortly that when these two subspaces coincide, we can conclude that $T$ is in fact, continuous.

\begin{cor}
$\,G(T)\,$ is a $\,\sigma$-closed subspace of $Y$.
\end{cor}

\begin{pro}
Suppose   $\,T:(X,\rho_{j})\rightarrow (Y,\sigma_{i})\,$ is linear and continuous.  Then  $\,G_{ji}(T)=Z_{\sigma_{i}}\,$.
\end{pro}
\noindent {\it Proof:}    Proposition 2.1-4 tells us that  $\,Z_{\sigma_{i}}<G_{ji}(T)$ always holds.  The continuity of $T$ implies the existence of $c>0$ such that$\,\sigma_{i}(T(x))\leq c \rho_{j}(x)\,$\, $\,\forall\, x\in X\,$.

\noindent  Let  $\,y\in G_{ji}(T)\,$ and let  $\,\{x_{n}\}_{\mathbb{N}}\stackrel{\rho_{j}}{\rightarrow} \hat{0}\in X\,$ such that  $\,\{T(x_{n})\}_{\mathbb{N}}\stackrel{\sigma_{i}}{\rightarrow} y\in Y\,$.   $\,\{T(x_{n})\}_{\mathbb{N}}\stackrel{\sigma_{i}}{\rightarrow} \tilde{0}\,$.\, These two convergences imply that $\,y=y-\tilde{0}\in Z_{\sigma_{i}}\,$, which proves that $\,G_{ji}(T) <Z_{\sigma_{i}}\,$.   \/ \/ $\Box$  \\

\begin{pro}
Let $\,T:(X,\rho_{j})\rightarrow (Y,\sigma_{i})\,$ be linear and such that  $\,G_{ji}(T)=Z_{\sigma_{i}}\,$ with  the spaces $\,(X,\rho_{j})\,$\ and  $\,(Y,\sigma_{i})\,$\, being complete.  Then $\,T\,$ is continuous.  
\end{pro}
\noindent {\it Proof} 
   We apply the closed graph theorem to the corresponding complete seminormed spaces.    Let  $\,\{x_{n}\}_{\Bbb{N}}\stackrel{\rho_{j}}{\rightarrow} x\in X\,$ and $\,\{T(x_{n})\}_{\Bbb{N}}\stackrel{\sigma_{i}}{\rightarrow} y\in Y\,$.   Then  $\,\{x_{n}-x\}_{\Bbb{N}}\stackrel{\rho_{j}}{\rightarrow} \hat{0}\in X\,$ and  $\,\{T(x_{n}-x)\}_{\Bbb{N}}\stackrel{\sigma_{i}}{\rightarrow} y-T(x)\in Y\,$.  That is, $\,y-T(x)\in G_{ji}(T)=Z_{\sigma_{i}}\,$, from which we obtain  $\,y\in T(x)+Z_{\sigma_{i}}\,$, and this tells us that the graph of $T$ is closed in $\,(X,\rho_{j})\times (Y,\sigma_{i})\,$.    By the closed graph theorem,$\,T:(X,\rho_{j})\rightarrow (Y,\sigma_{i})\,$ is continuous. \/ \/  $\Box$  \\

\begin{cor} \begin{enumerate}
\item  If $T:(X,\tau)\rightarrow (Y,\sigma)$ is linear and continuous, then for each $\,\sigma_{i}\,$ there is a $\,\rho_{j}\,$\ such that  $\,G_{ji}(T)=Z_{\sigma_{i}}\,$.
\item   If  for each  $\sigma_{i}$ the seminormed space $\,(Y,\sigma_{i})\,$ is complete and if  there exists $\,\rho_{j}\,$ such that  $\,G_{ji}(T)=Z_{\sigma_{i}}\,$\,  and $\,(X,\rho_{j})\,$ is also complete, then \, $\,T\,$\, is continuous.
\end{enumerate}
\end{cor}

\noindent {\it Proof} 
  1.  The continuity of  $\,T\,$\, implies that for each seminorm, $\sigma_{i}$ from the family that generates $ \sigma\,$\, there is a corresponding seminorm $\,\rho_{j}\in \tau$ such that $\,T:(X,\rho_{j})\rightarrow (Y,\sigma_{i})\,$\, is continuous. We apply  Proposition 2.3  to obtain the conclusion.  For item 2, we directly apply  Proposition  2.3  \/ \/  $\Box$  \\ 

\begin{cor}  
For any linear $\,T: (X,\tau) \rightarrow (Y,\sigma)\,$, with $(X,\tau)$ and $(Y,\sigma)$ complete,   one has  $G(T)=\{\tilde{0}\}\,$

\end{cor}
\noindent {\it Proof} 
   $\{\tilde{0}\}<G(T)=\bigcap_{J\times I}G_{ji}(T)<\bigcap_{i\in I}Z_{\sigma_{i}}=\{\tilde{0}\}\,$.  \/ \/  $\Box$   \\  
\bigskip

\noindent Observe that if $T:(X,\rho_{j_{o}})\rightarrow (Y,\sigma_{i_{o}})$ is continuous,  Proposition 2.3  tells us that  $\,G_{j_{o}i_{o}}=Z_{\sigma_{i_{o}}}$. On the other hand, 

$$\,Z_{\sigma_{i_{o}}}=G_{j_{o}i_{o}}(T)=\bigcap_{j\in J}G_{ji_{o}}(T)=\bigcap_{i\in I}G_{j_{o}i}(T)<G_{j_{o}i_{o}}(T).  \; \; \Box$$

We conclude that the above subspaces are equal.  

\section{A Continuity result.}

\noindent  The main result of this section is that, under our constructions,  we can conclude the continuity of  $T:(X,\tau)\rightarrow (Y,\sigma)\,$  whenever the closed graph theorem applies. This generalizes Sinclair  \cite{SC} to large classes of  locally convex spaces, such as Pt\'{a}k spaces, webbed spaces, and so on. See  \cite{PC}.    \\
\begin{teo}
If $G(T)=\{\tilde{0}\}$ and the spaces  $(X,\rho_{j})\,$ and  $\,(Y,\sigma_{i})\,$ are complete for each  $\,(j,i)\in J\times I\,$, then the graph of $T:(X,\tau)\rightarrow (Y,\sigma)\,$ is   $\,(\tau \times \sigma)$- closed in\, $\,X\times Y$.  In this case, a linear $T:(X,\tau)\rightarrow (Y,\sigma)$ is continuous whenever the lcs  $(X,\tau)\,$ and  $\,(Y,\sigma)$ satisfy the  assumptions of a closed graph theorem.  
\end{teo}

\noindent {\it Proof} 
  Let $\Delta$ be a directed set and let $\,\{(x_{\alpha},T(x_{\alpha}))\}_{\Delta}\stackrel{\tau \times \sigma}{\rightarrow} (x,y)\in X\times Y\,$.   That is, $\,\{x_{\alpha}\}_{\Delta}\stackrel{\tau}{\rightarrow} x\,$ and  $\,\{T(x_{\alpha})\}_{\Delta}\stackrel{\sigma}{\rightarrow} y\,$.  Let $(j,i)\in J\times I\,$ be fixed, but arbitrary.  Given that the topologies induced by the seminorms  $\,\rho_{j}\,$\, y\, $\,\sigma_{i}\,$ in $\,X$ and $Y$, respectively, are weaker than $\,\tau\,$\, y\, $\,\sigma\,$ (respectively), we have that $\,\{x_{\alpha}\}_{\Delta}\stackrel{\rho_{j}}{\rightarrow} x$ and $\{T(x_{\alpha})\}_{\Delta}\stackrel{\sigma_{i}}{\rightarrow} y$.  \\

\noindent   Nevertheless, we need convergence in terms of sequences.  For each  $\,n\in \mathbb{N}$, put $\epsilon_{n}=\frac{1}{n}$.  Thus, there exist $\,\alpha_{n}'\,$ and $\,\alpha_{n}''\,$ in  $\,\Delta\,$ such that $(\forall\, \alpha \geq \alpha_{n}') \, \rho_{j}(x_{\alpha}-x)<\epsilon_{n}$, and $(\forall\, \alpha \geq \alpha_{n}'' \sigma_{i}) \;  T((x_{\alpha})- y) <\epsilon_{n}\,$.   \
 
\noindent Let $\,\alpha_{1}\in \Delta$ such that $\,\alpha_{1}\geq \alpha_{1}',\,\alpha_{1}''\,$.    Then let $\,\alpha_{2}\in \Delta\,$ such that  $\,\alpha_{2}\geq \alpha_{2}',\,\alpha_{2}'',\, \alpha_{1}\,$.  In general, we have   $\alpha_{n}\in \Delta\,$ such that  $\,\alpha_{n}\geq \alpha_{n}',\,\alpha_{n}'',\, \alpha_{n-1}\,$.

\noindent   We now have two sequences in place of subnets, namely,   $\,\{x_{\alpha_{n}}\}_{\Bbb{N}}\stackrel{\rho_{j}}{\rightarrow} x\in X\,$ and  $\,\{T(x_{\alpha_{n}})\}_{\Bbb{N}}\stackrel{\sigma_{i}}{\rightarrow} y\in Y\,$.   \

\noindent Following the proof of  Proposition 2.3  we obtain  $\,y-T(x)\in G_{ji}(T)\,$, though not necessarily in   $\,Z_{\sigma_{i}}\,$.  By the above,  $\,y-T(x)\in \bigcap_{J\times I}G_{ji}(T)=G(T)=\{\tilde{0}\}\,$, which implies $\,y=T(x)$, and therefore, $\,(x,y)\,$ belongs to the graph of $\,T\,$ in $\,X\times Y\,$.    \/ \/  $\Box$   \\

\noindent  The following is an example of how this result applies to generalize a result  like  \cite{SC}, Lemma 1.3 - (i), page  8:  

\begin{cor}
Suppose  $\,G(T)=\{\tilde{0}\}\,$ and the lcs  $\,(X,\rho_{j})\,$\, $\,(Y,\sigma_{i}))\,$\,  are  complete for each  $\,(j,i)\in J\times I\,$. If the spaces  $\,(X,\tau)\,$,\, $\,(Y,\sigma)\,$ are, respectively,  barrelled and  Pt{\'a}k spaces, then the linear map  $\,T:(X,\tau)\rightarrow (Y,\sigma)\,$\, is continuous.\, 

\end{cor}

\noindent {\it Proof} 
     See  \cite{NB}  14.9.1, page 327.     \/ \/  $\Box$   

\section{Eventually constant intertwining linear maps.}  

\noindent  Our goal in this section is to generalize Sinclair's fundamental  result  \cite{SC},  Lemma 1.6, regarding eventually constant intertwining linear maps between  Banach spaces, to linear maps between  complete, locally convex spaces.  We need some preliminary results as we build  from the ground up using seminorms.  

\begin{pro}
Let $\,R:(Y,\sigma_{i})\rightarrow (W,\lambda_{k})\,$, be linear and continuous. If there is $j\in J$  such that  $\,R(G_{ji}(T))<Z_{\lambda_{k}}$, then 
$\,R\circ T:(X,\rho_{j})\rightarrow (W,\lambda_{k})\,$\, is continuous. 

\end{pro}

\noindent {\bf Comment}:   The requirement that $\,R(G_{ji}(T))<Z_{\lambda_{k}}\,$ would make for a trivial result, since it would force the null sequences in the domain to satisfy that their images under the composition $\,R \circ T\,$, to converge to the zero in the codomain, which of course, would be continuity.  The fact that  $\,\{T(x_{n})\}_{\Bbb{N}}$ does not always converge in $Y$ is what makes this result interesting, though it also makes for more work to prove it.  \\   

\noindent {\it Proof} 
   We will use   Proposition 2.3, proving that  $\,G_{jk}(R\circ T)=Z_{\lambda_{k}}$.  By part 4 of  Proposition 2.1,  $Z_{\lambda_{k}}<G_{jk}(R\circ T)\,$.   \

Let $\,w\in G_{jk}(R\circ T)<W\,$\ and let   $\,\{x_{n}\}_{\Bbb{N}}\stackrel{\rho_{j}}{\rightarrow} \hat{0}\in X\,$ such that  $\,\{(R\circ T)(x_{n})\}_{\Bbb{N}}\stackrel{\lambda_{k}}{\rightarrow} w\in Y\,$.   If $\,\{T(x_{n})\}_{\Bbb{N}}\stackrel{\sigma_{i}}{\rightarrow} \tilde{0}\,$, then $\,w\in G_{ik}=Z_{\lambda_{k}}\,$.  \

If $\,\{T(x_{n})\}_{\Bbb{N}}\stackrel{\sigma_{i}}{\rightarrow} y\neq \tilde{0}\,$, then  $\,y\in G_{ji}(T)\,$, from which we have $\,R(y)\in R(G_{ji}(T))\in Z_{\lambda_{k}}\,$.  Thus,  $\,\{R(T(x_{n})-y)\}_{\Bbb{N}}\stackrel{\lambda_{k}}{\rightarrow} w-R(y)$.  From this, $\,w-R(y)\in G_{ik}(R)=Z_{\lambda_{k}}\,$ and  $\,w\in R(y)+Z_{\lambda_{k}}=Z_{\lambda_{k}}\,$.   \
Let  $\,\eta:(Y,\sigma_{i})\rightarrow (Y/G_{ji}(T),\hat{\sigma}_{i})\,$ be the canonical epimorphism. \
 Let   $\,\hat{R}:(Y/G_{ji}(T),\hat{\sigma}_{i})\rightarrow (W,\lambda_{k})\,$ under the association  $\,\hat{R}(\eta(y))=R(y)\,$.  This function is well defined, and by hypothesis, $\,R(G_{ji}(T))<Z_{\lambda_{k}}=G_{ik}(R)\,$.    \\

\noindent {\bf Claim 1}:\, $\,\hat{R}\,$\, is continuous; that is,  $\,G_{ik}(\hat{R})=Z_{\lambda_{k}}\,$.

\noindent  Let  $\,w\in G_{ik}(\hat{R})\,$ and let   $\,\{\eta(y_{n})\}_{\Bbb{N}}\stackrel{\hat{\sigma}_{i}}{\rightarrow} \eta(\tilde{0})\in Y/G_{ji}(T)\,$ such that  $\,\{R(y_{n})\}_{\Bbb{N}}=\{\hat{R}(\eta(y_{n}))\}_{\Bbb{N}}\stackrel{\lambda_{k}}{\rightarrow} w\in W\,$.  Then there exists $\,\{y_{n}'\}_{\Bbb{N}}\subset Y\,$ such that  $\,\eta(y_{n}')=\eta(y_{n})\,$\, y\, $\,\{y_{n}'\}_{\Bbb{N}}\stackrel{\sigma_{i}}{\rightarrow} \tilde{0}\,$.

\noindent    From this, we have  $\,\{\hat{R}(\eta(y_{n}))\}_{\Bbb{N}}=\{\hat{R}(\eta(y_{n}'))\}_{\Bbb{N}}=\{R(y_{n}')\}_{\Bbb{N}}\stackrel{\lambda_{k}}{\rightarrow} 0\in W\,$.  These two convergences imply that  $\,w\in Z_{\lambda_{k}}\,$.    \\
 
\noindent  {\bf Claim 2}:\, $\,\eta \circ T : (X,\rho_{j})\rightarrow (Y/G_{ji}(T),\hat{\sigma}_{i})\,$\, is also continuous;   that is,   $\,G_{ji}(\eta \circ T)=Z_{\lambda_{k}}\,$.

\noindent Let  $\,\eta(y)\in G_{ji}(\eta \circ T)\,$ and let  $\,\{x_{n}\}_{\Bbb{N}}\stackrel{\rho_{j}}{\rightarrow} \hat{0}\in X\,$ such that  $\,\{(\eta \circ T)(x_{n})\}_{\Bbb{N}}\stackrel{\hat{\sigma}_{i}}{\rightarrow} \eta(y)\in Y/G_{ji}(T)\,$.  We now have  $\,\{\eta(T(x_{n})-y)\}_{\Bbb{N}}\rightarrow \eta(\tilde{0})\,$, in other words,  $\,\hat{\sigma}_{i}(\eta(T(x_{n})-y))\rightarrow 0\,$.  \
 Thus, for each  $\,n\in \Bbb{N}\,$\, there is \, $\,y_{n}\in G_{ji}(T)\,$ such that  $\,\hat{\sigma}_{i}(\eta(T(x_{n})-y))\leq \sigma_{i}(T(x_{n})-y-y_{n})< \hat{\sigma}_{i}(\eta(T(x_{n})-y))+\frac{1}{n}\,$. 

\noindent  $\,y_{n}\in G_{ji}(T)\, \Rightarrow\, \eta(T(x_{n})-y_{n})=\eta(T(x_{n}))\,$.  This means that, given $\,\epsilon >0\,$,\, there exists \, $\,N\in \Bbb{N}\,$ for which  $\,\forall\, n\geq N\,$: 
$$\,\hat{\sigma}_{i}(\eta(T(x_{n})-y)) \leq \sigma_{i}(T(x_{n})-y-y_{n})  <\hat{\sigma}_{i}(\eta(T(x_{n})-y))+ \frac{1}{n}< \frac{\epsilon}{2}+\frac{\epsilon}{2}.    $$

\noindent From this,  $\,\{T(x_{n})-y_{n}\}_{\Bbb{N}}\stackrel{\sigma_{i}}{\rightarrow} y\,$.  On the other hand,  for each $y_{n}\in G_{ji}(T)\,$  there exists   $\,\{x_{nm}\}_{m\in \Bbb{N}}\stackrel{\rho_{j}}{\rightarrow} \hat{0}\in X\,$, as $\,m\rightarrow \infty\,$,  with  $\,\{T(x_{nm})\}_{m\in \Bbb{N}}\stackrel{\sigma_{i}}{\rightarrow} y_{n}\,$.  For each  $\,n\in \Bbb{N}\,$\, there exists \, $\,m_{n}\,$ such that $\,\rho_{j}(x_{nm})<\frac{1}{n}$ and 

$$  \sigma_{i}(T(x_{nm})-y_{n})<\frac{1}{n}, \; \;    \forall\, m\geq m_{n}.   $$

\noindent  We thus have,  $\,\{x_{n}-x_{nm_{n}}\}_{n\in \Bbb{N}}\rightarrow \hat{0}\,$ and 

\begin{eqnarray}   \sigma_{i}(T(x_{n}-x_{nm_{n}})-y)& = &\sigma_{i}(T(x_{n})-T(x_{nm_{n}})-y) \nonumber  \\
                            & \leq &  \sigma_{i}(T(x_{n})-y_{n}-y)+\sigma_{i}(y_{n}-T(x_{nm_{m}}))  \nonumber  \\     
                             &  \rightarrow & 0.   \nonumber  \end{eqnarray}  

\noindent The above tells us that $\,y\in G_{ji}(T)\,$, and therefore,  $\,\eta(y)=\eta(\tilde{0})=Z_{\lambda_{k}}\,$.    \\

\noindent      Having shown that the functions  $\,\eta\circ T\,$ and $\,\hat{R}\,$ are continuous, their composition   $\,\hat{R}\circ \eta\circ T\,$ is of course,  also continuous.  

\noindent     Because $\,\hat{R}\circ \eta \circ T= R\circ T$, $R\circ T\,$\, is continuous, which is what we intended to prove.    \/ \/ $\Box$   

\begin{rem}
If the linear function $\,T:(X,\rho_{j})\rightarrow (Y,\sigma_{i})\,$ were continuous, the number $M_{ji}$  would exist, where   $\,0\leq M_{ji}= \inf\,\{\sigma_{i}(Tx)\,|\,\rho_{j}(x)\leq 1\}\,$ which satisfies $\,\sigma_{i}(T(x))\leq M_{ji} \rho_{j}(x)\,$ for each  $\,x\in X\,$. Analogously, if the linear function $\,R:(Y,\sigma_{i})\rightarrow (W,\lambda_{k})\,$ is continuous, then $\,0\leq M_{ik}= \inf\, \{\lambda_{k}(R(y))\,|\,\sigma_{i}(y)\leq 1\}\,$ exists such that   $\,\lambda_{k}(R(y))\leq M_{ik} \sigma_{i}(y)\,$ for each  $\,y\in Y\,$.\, We conclude from this that for the composition  $\,(R\circ T): (X,\rho_{j})\rightarrow (W,\lambda_{k})\,$ we have  $\,\lambda_{k}((R\circ T)(x))\leq M_{ik} M_{ji} \rho_{j}(x)\,$ for each $\,x\in X\,$.   
\end{rem}

\noindent     We can ask the question:  What happens if  $\,R\circ T\,$\, is continuous, but $T$ is not, as  in the context of Proposition  4.1?     We will see what  this tells us  in the next proposition, after the following technical lemma.      

\begin{lem}
\begin{eqnarray}  M_{ik} & =&  \inf\, \{\lambda_{k}(R(y))\,|\,\sigma_{i}(y)\leq 1\}   \nonumber \\
                                        &= & M_{\hat{i}k} \nonumber  \\
                                         &  =  &  \inf\,\{\lambda_{k}(R(y+G_{ji}(T))\,|\, \hat{\sigma}_{i}(y+G_{ji}(T))\leq 1\}. \nonumber
  \end{eqnarray}  
\end{lem}
\noindent {\it Proof} 
   From Proposition 4.1 we have that the function   $\,\hat{R} : (Y/G_{ji}(T), \hat{\sigma}_{i}) \rightarrow (W,\lambda_{k})\,$,  $\,\hat{R}(y+G_{ji}(T))=R(y)\,$, is continuous. Hence, there exists  a number   $\,M_{\hat{i}k} \geq 0$ for which we have  
$$ \lambda_{k}(R(y))= \lambda_{k}(\hat{R}(y+G_{ji}(T)))\leq M_{\hat{i}k}\hat{\sigma}_{i}(y+G_{ji}(T))\leq  M_{\hat{i}k}\sigma_{i}(y), $$     

\noindent     $\,\forall \,y\in Y\,$.  The previous statement tells us that $\,M_{ik}\leq M_{\hat{i}k}\,$, where   

$$  M_{ik}= \inf\,\{\lambda_{k}(R(y))\,|\,\sigma_{i}(y)\leq 1\}.  $$

\noindent   We recall that $\,\hat{\sigma}_{i}(y+G_{ji}(T))= \inf\,\{\sigma_{i}(y+z)\,|\,z\in G_{ji}(T)\}\leq \sigma_{i}(y)\,$.  By supposing  $\,M_{ik}< M_{\hat{i}k}\,$ we have that there exists  $\,y_{1}\in Y\,$ with  $\,\sigma_{i}(y_{1})\leq 1\,$ such that   $\,M_{ik} \leq \lambda_{k}(R(y_{1}))< M_{\hat{i}k}\,$.   \\

\noindent     Nevertheless, $\,\lambda_{k}(R(y_{1}))=\lambda_{k}(\hat{R}(y_{1}+G_{ji}(T))\,$ and  $\,\hat{\sigma}_{i}(y_{1}+G_{ji}(T))\leq \sigma_{i}(y_{1})\leq 1\,$, that is,  $\,M_{\hat{i}k}\leq \lambda_{k}(R(y_{1}))< M_{\hat{i}k}\,$.  This contradiction shows that $\,M_{ik}=M_{\hat{i}k}\,$.   \/ \/ $\Box$  \\  

\noindent  Here is the proposition alluded to a short time ago:  
\begin{pro}
Suppose the linear function $\,R:(Y,\sigma_{i})\rightarrow (W,\lambda_{k})\,$\,  is continuous,  such that  $\,(R\circ T): (X,\rho_{j})\rightarrow (W,\lambda_{k})\,$\,is also continuous.  Suppose additionally, that  $\,R(G_{ji}(T))<Z_{\lambda_{k}}\,$.  Then there exists a number     $M>0$, independent of the the function $R$ and of the vector space $W$, such that    $\,\lambda_{k}((R\circ T)(x))\leq M_{ik} M \rho_{j}(x)\,$ for each   $x\in X$.

\end{pro}
\noindent {\it Proof} 
  Based on what we did in Proposition  4.1, we have that  $\,R \circ T=\hat{R} \circ \eta \circ T\,$, where  $\,\eta \circ T :(X,\rho_{j})\rightarrow (Y/G_{ji}(T),\hat{\sigma}_{i})\,$\,   is a linear, continuous function, as well as  $\,\hat{R}: (Y/G_{ji}(T),\hat{\sigma}_{i}) \rightarrow (W,\lambda_{k})\,$.   By the continuity of  $\,\eta \circ T\,$,\, there is a number    $\,M_{j \hat{i}}\geq 0\,$ such that 

 $$\,\hat{\sigma}_{i}((\eta \circ T)(x))\leq M_{j\hat{i}}\rho_{j}(x) \; \;  (\forall \, x\in X).   $$

\noindent   Finally, $(\forall \,x\in X)$, 
\begin{eqnarray}  \lambda_{k}((R\circ T)(x))& =& \lambda_{k}((\hat{R}\circ \eta \circ T)(x))  \nonumber \\
                                                                     &   \leq & M_{ik}\hat{\sigma}_{i}((\eta \circ T)(x))  \nonumber  \\
                                                                     & \leq & M_{ik}M_{j\hat{i}}\rho_{j}(x),  \nonumber  
\end{eqnarray} 

\noindent     Choose  $\,M=M_{j\hat{i}}\,$\, and observe that this number is independent of the function $R$ as well as the vector space $W$.   \/ \/ $\Box$

\begin{cor}
Suppose the function   $\,R:(Y,\sigma)\rightarrow (W,\lambda)\,$\, is continuous, such that  $\,(R\circ T): (X,\tau)\rightarrow (W,\lambda)\,$\,is also  continuous.  Moreover, because for each      $\,\lambda_{k}\in \lambda\,$\, there exists\, $\,\sigma_{i}\in \sigma\,$ for which  $\,R:(Y,\sigma_{i})\rightarrow (W,\lambda_{k})\,$\, is continuous, we can have that for this    $\,\sigma_{i}\in \sigma\,$ there will  exist $\,\rho_{j}\in \tau\,$ such that  $\,R(G_{ji}(T))<Z_{\lambda_{k}}\,$. Then,\, $\,(R \circ T):(X,\tau)\rightarrow (W,\lambda)\,$\, is continuous.   
\end{cor}

\vspace{3mm}
\noindent  The next five results inform us about the vector subspaces of     $W$ :\, $\,G_{jk}(R\circ t)\,$,\, $\,G_{ik}(R)\,$ and  $\,R(G_{ji}(T))\,$.   \\

\begin{cor}
Let $\,R:(Y,\sigma)\rightarrow (W,\lambda)\,$\,  be linear and continuous,  that is, for each $\,\lambda_{k}\in \lambda\,$\,there exists $\,\sigma_{i}\in \sigma\,$ such that   $\,R:(Y,\sigma_{i})\rightarrow (W,\lambda_{k})\,$\, is continuous.  If for this seminorm  $\,\sigma_{i}\,$\, there is $\,j\in J\,$ for which  $\,R(G_{ji}(T))<Z_{\lambda_{k}}\,$,\, then   \, $\,R\circ T:(X,\tau)\rightarrow (W,\lambda)\,$\, is continuous.    
\end{cor}

\begin{pro}
Let $\,R:(Y,\sigma_{i})\rightarrow (W,\lambda_{k})\,$\, be linear and continuous. Then  $\,R(G_{ji}(T))<G_{ik}(R)=Z_{\lambda_{k}}\,$ and  $\,\overline{R(G_{ji}(T))}^{\lambda_{k}}<G_{ik}(R\circ T)\,$ for every  $\,j\in J\,$.
\end{pro}

\noindent {\it Proof} 
    Consider  $\,j\in J\,$ and $\,T:(X,\rho_{j})\rightarrow (Y,\sigma_{i})\,$,\, linear but not necessarily continuous.  By the continuity of $R$, Proposition 2.3 tells us that  $\,G_{ik}(R)=Z_{\lambda_{k}}\,$. Let $y\in G_{ji}(T)\,$ and let   $\,\{x_{n}\}_{\Bbb{N}}\stackrel{\rho_{j}}{\rightarrow} \hat{0}\in X\,$ such that  $\,\{T(x_{n})\}_{\Bbb{N}}\stackrel{\sigma_{i}}{\rightarrow} y\in Y\,$.\, Then,\, $\,\{R(T(x_{n}))\}_{\Bbb{N}}\stackrel{\lambda_{k}}{\rightarrow} R(y)\,$ and $\,(R\circ T)(x_{n})\}_{\Bbb{N}}\stackrel{\lambda_{k}}{\rightarrow} 0\in W\,$.  These two convergences imply two things.  The first is that $\,R(y)\in Z_{\lambda_{k}}=G_{ik}(R)\,$ and the other is $\,R(y)\in G_{jk}(R\circ T)\,$. From the first one we conclude that  $\,R(G_{ji}(T))<G_{ik}(R)\,$, and the second tells us that  $\,R(G_{ji}(T))<G_{jk}(R\circ T)\,$.  From this, 

$$   \overline{R(G_{ji}(T))}^{\lambda_{k}}< \overline{G_{jk}(R\circ T)}^{\lambda_{k}}= G{jk}(R\circ T).  \; \; \Box $$  


\begin{cor}
Suppose  $\,R:(Y,\sigma_{i})\rightarrow (W,\lambda_{k})\,$ is linear and continuous, and  $\,T:(X,\rho_{j})\rightarrow (Y,\sigma_{i})\,$ is linear, such that the composition   $(R\circ T):(X,\rho_{j})\rightarrow (W,\lambda_{k})\,$ is continuous.   Then  $G_{ik}(R)=Z_{\lambda_{k}}=G_{jk}(R\circ T)\,$.
\end{cor}
 
\begin{cor}
Let $\,R:(Y,\sigma)\rightarrow (W,\lambda)\,$ be linear and continuous.   If  $\,R\circ T:(X,\tau)\rightarrow (W,\lambda)\,$\, is continuous, then for each  $\,\lambda_{k}\in \lambda\,$\, there exists \, $\,j\in J\,$ and  $\,i\in I\,$\ for which  $\,\overline{R(G_{ji}(T))}^{\lambda_{k}}<G_{ik}(R)=Z_{\lambda_{k}}=G_{jk}(R\circ T)\,$.
\end{cor}

\begin{pro}
Let  $\,R:(Y,\sigma_{i})\rightarrow (W,\lambda_{k})\,$ be linear and continuous.  Then  

$$\,\overline{R(G_{ji}(T))}^{\lambda_{k}}=G_{jk}(R\circ T).  $$

\end{pro}
\noindent {\it Proof} 
   Let  $\,y\in G_{ji}(T)\,$ and let  $\,\{x_{n}\}_{\Bbb{N}}\stackrel{\rho_{j}}{\rightarrow} \hat{0}\in X\,$ such that 

$$\,\{T(x_{n})\}_{\Bbb{N}}\stackrel{\sigma_{i}}{\rightarrow} y\in Y.  $$

\noindent  By the continuity of $R$,  $\,\{R(T(x_{n}))\}_{\Bbb{N}}\stackrel{\lambda_{k}}{\rightarrow} R(y)\,$, from which   $\,R(y)\in G_{jk}(R\circ T)\,$.  In other words,   $\,R(G_{ji}(T))<G_{jk}(R\circ T)\,$ and because  $\,G_{jk}(R\circ T)\,$ is $\,\lambda_{k}$- closed, we have 
  
$$  \overline{R(G_{ji}(T))}^{\lambda_{k}}<G_{jk}(R\circ T).  $$

\noindent    By Proposition 2.4,   $\,R\circ T:(X,\rho_{j})\rightarrow (W,\lambda_{k})\,$\, is continuous.  \\

\noindent  Now we prove the other contention:  By the continuity of   $\,R\circ T\,$, we have $\,Z_{\lambda_{k}}=G_{jk}(R\circ T)\,$. Consider the canonical epimorphism  $\,\eta:W\rightarrow (W/\overline{R(G_{ji}(T))},\hat{\lambda}_{k})\,$, which happens to be continuous.  That continuity implies that     $\, \eta \circ R:(Y,\sigma_{i})\rightarrow  (W/\overline{R(G_{ji}(T))},\hat{\lambda}_{k})\,$, is continuous and that also, 

 $$ (\eta \circ R)\circ T=\eta \circ (R\circ T):(X,\rho_{j})\rightarrow  (W/\overline{R(G_{ji}(T))},\hat{\lambda}_{k}).  $$

\noindent  By  Proposition 4.1  $\,\overline{R(G_{ji}(T))}^{\lambda_{k}}< G_{ik}(R)=Z_{\lambda_{k}}\,$.  Let  $\,y\in G_{jk}(R\circ T)=Z_{\lambda_{k}}\,$ and  $\,\{x_{n}\}_{\Bbb{N}}\stackrel{\rho_{j}}{\rightarrow} \hat{0}\,$ such that  $\,\{(R\circ T)(x_{n})\}_{\Bbb{N}}\stackrel{\lambda_{k}}{\rightarrow} y\,$.  We then have two convergences, knowing that  $\,\{(\eta\circ R\circ T)(x_{n})\}_{\Bbb{N}}\rightarrow \eta(y)\,$ and that  $\,\{(\eta\circ R\circ T)(x_{n})\}_{\Bbb{N}}\rightarrow \eta(\hat{0})\,$, from which  $\,\eta(y)\in Z_{\hat{\lambda}_{k}}=\eta(0)=[0]\,$. This implies that  $\,y\in \overline{R(G_{ji}(T))}^{\lambda_{k}}\,$. \/ \/ $\Box$    

\begin{cor}
Let $\,R:(Y,\sigma_{i})\rightarrow (W,\lambda_{k})\,$ be linear and continuous, and $\,T:(X,\rho_{j})\rightarrow (Y,\sigma_{i})\,$ linear such that the composition  $\,(R\circ T):(X,\rho_{j})\rightarrow (W,\lambda_{k})\,$ is continuous.  Then   $\,\overline{R(G_{ji}(T)}^{\lambda_{k}}=Z_{\lambda_{k}}\,$.


\end{cor}
\noindent {\it Proof} 
   By Proposition 4.10  $\,\overline{R(G_{ji}(T)}^{\lambda_{k}}=G_{jk}(R\circ T)=Z_{\lambda_{k}}\,$, this last equality being from the continuity of $\,R\circ T\,$. \/ \/ $\Box$   

\begin{cor}
Let $\,R:(Y,\sigma)\rightarrow (W,\lambda)\,$ be linear and continuous.  Then
$\,\overline{R(G(T))}^{\lambda}<G(R\circ T)\,$.  This is  Proposition 4.4, with the topologies  $\,\tau\,$ and  $\,\sigma\,$.
\end{cor}
\noindent {\it Proof} 
\begin{eqnarray}  \overline{R(G(T))}^{\lambda}& =& \overline{R(\bigcap_{J\times I}G_{ji}(T))}^{\lambda} \nonumber \\
                                                                             & < & \overline{\bigcap_{J\times I}R(G_{ji}(T))}^{\lambda}  \nonumber \\
                                                                            & < &  \bigcap_{J\times I}\overline{R(G_{ji}(T))}^{\lambda}  \nonumber \\
                                                                             &=&  \bigcap_{J\times K}G_{jk}(T)=G(R\circ T) \nonumber   \;   \; \Box
\end{eqnarray} 

\begin{pro}
Suppose  $\,S:(X,\rho_{l})\rightarrow (X,\rho_{j})\,$ and  $\,R:(Y,\sigma_{i})\rightarrow (Y,\sigma_{k})\,$\, are continuous, linear functions, and consider  $\,T:(X,\rho_{j})\rightarrow (Y,\sigma_{i})$.  These satisfy the following:   
\begin{enumerate}

\item  $\,G_{li}(T\circ S)< G_{ji}(T)\,$ 
\item  $\,R(G_{ji}(T\circ S))< G_{jk}(R\circ T)\,$ 
\item If  $\,T\circ S=R\circ T\,$, then $\,R(G_{ji}(T))< G_{jk}(T)\,$.   
\end{enumerate}
\end{pro}
\noindent {\it Proof} 
  1) Let  $\,y\in G_{li}(T\circ S)\,$, let   $\,\{x_{n}\}_{\Bbb{N}}\stackrel{\rho_{j}}{\rightarrow} \hat{0}\in X\,$, and $\,\{(T\circ S)(x_{n})\}_{\Bbb{N}}\stackrel{\sigma_{i}}{\rightarrow} y\in Y\,$.  By the continuity of  $S$,\, $\,\{S(x_{n})\}_{\Bbb{N}}\stackrel{\rho_{j}}{\rightarrow} \hat{0}\,$, $\,\{T(S(x_{n}))\}_{\Bbb{N}}\stackrel{\sigma_{i}}{\rightarrow} y\,$.   This proves that  $\,y\in G_{ji}(T)\,$, which in turn implies  $\,G_{li}(T\circ S)<G_{ji}(T)\,$.   \

\noindent  2)  Let  $\,y\in G_{ji}(T\circ S)\,$.  From this,  there exists    $\,\{x_{n}\}_{\Bbb{N}}\stackrel{\rho_{l}}{\rightarrow} \hat{0}\in X\,$ such that   $\,\{(T\circ S)(x_{n})\}_{\Bbb{N}}\stackrel{\sigma_{i}}{\rightarrow} y\in Y\,$.    By the continuity of  $\,S\,$,\, $\,\{S(x_{n})\}_{\Bbb{N}}\stackrel{\rho_{j}}{\rightarrow} \hat{0}\,$.  By the continuity of   $\,R\,$,\, $\,\{(R\circ T)(S(x_{n}))\}=\{R(T\circ S)(x_{n})\}_{\Bbb{N}}\stackrel{\sigma_{k}}{\rightarrow} R(y)\,$.    From here, we obtain  $\,R(y)\in G_{jk}(R\circ T)\,$.    \

\noindent   3)   Let    $\,y\in G_{ji}(T)\,$ and let $\,\{x_{n}\}_{\Bbb{N}}\stackrel{\rho_{j}}{\rightarrow} \hat{0}\in X\,$\, y\, $\,\{T(x_{n})\}_{\Bbb{N}}\stackrel{\sigma_{i}}{\rightarrow} y\in Y\,$.  This gives us  $\,\{S(x_{n})\}_{\Bbb{N}}\stackrel{\rho_{j}}{\rightarrow} \hat{0}\,$ and  $\,\{T(S(x_{n}))\}_{\Bbb{N}}=\{R(T(x_{n}))\}_{\Bbb{N}}\stackrel{\sigma_{k}}{\rightarrow} R(y)\,$.\,   The last piece tells us that  $\,R(y)\in G_{jk}(T)\,$. \/ $\Box$   
\bigskip

\begin{cor}
If  $\,S:(X,\tau)\rightarrow (X,\tau)\,$ and  $\,R:(Y,\sigma)\rightarrow (Y,\sigma)\,$\,  are continuous, linear functions, then they satisfy the following:   
\begin{enumerate}
\item    $\,G(T\circ S)< G(T)\,$ 
\item   $\,R(G(T\circ S))< G(R\circ T)\,$ 
\item  If $\,T\circ S=R\circ T\,$,\, then \, $\,R(G(T))< G(T)\,$. 
\end{enumerate}
\end{cor}

\begin{pro}
Let   $\,U<X\,$ and  $\,V<Y\,$\, be closed vector subspaces for which     $\,T(U)<V\,$.  Let  $\,\hat{T}:(X/U,\hat{\rho}_{j}) \rightarrow (Y/V,\hat{\sigma}_{i})\,$ be linear, such that $\,\hat{T}\circ \eta_{U}= \eta_{V}\circ T\,$. (This is possible by the hypothesis).   Then $\,\hat{T}\,$\, is continuous iff $\,G_{ji}(T)<V\,$. \\ 
$\,\eta_{U}:(X,\rho_{j})\rightarrow (X/U,\hat{\rho}_{j})\,$ and  $\,\eta_{V}:(Y,\sigma_{i})\rightarrow (Y/V,\hat{\sigma}_{i})\,$ are the canonical epimorphisms, which are continuous. 
\end{pro}
  ($ \Rightarrow$ ): \/ The continuity of $\,\hat{T}\,$ implies that  $\,\hat{T}\circ \eta_{U}\,$\, is also continuous.   Let  $\,y\in G_{ji}(T)\,$ as well as $\,\{x_{n}\}_{\Bbb{N}}\stackrel{\rho_{j}}{\rightarrow} \hat{0}\in X\,$, such that  $\,\{T(x_{n})\}_{\Bbb{N}}\stackrel{\sigma_{i}}{\rightarrow} y\in Y\,$.  This shows that  $\{\eta_{U}(x_{n})\}_{\Bbb{N}}\rightarrow \eta_{U}(\hat{0})\,$  and  $\,\{\hat{T}(\eta_{U}(x_{n})\}_{\Bbb{N}}=\{\eta_{V}(T(x_{n})\}_{\Bbb{N}}\rightarrow \eta_{V}(y)\,$. Meanwhile,  $\,\{\hat{T}(\eta_{U}(x_{n})\}_{\Bbb{N}}\rightarrow \hat{T}(\eta_{U}(\hat{0})=\eta_{V}(T(\hat{0}))=\eta_{V}(\tilde{0})\,$.  These two convergences tell us that $\,y\in V\,$.    \\

\noindent ($ \Leftarrow$): \/ Our task is to show that  $\,G_{ji}(\hat{T})=Z_{\hat{\sigma}_{i}}=\eta_{V}(\tilde{0})\,$.  To this end, we start with  $\,G_{ji}(T)<V\,\Rightarrow\, \eta_{V}(G_{ji}(T))=\eta_{V}(\tilde{0})<Z_{\hat{\sigma}_{i}}\,$.  By  Proposition 4.1,\, $\,\eta_{V}\circ T\,$\, is continuous, and thus, so is  $\,\hat{T}\circ \eta_{U}\,$.   Let  $\,\eta_{V}(y)\in G_{ji}(\hat{T})\,$  and let  $\,\{\eta_{U}(x_{n})\}_{\Bbb{N}}\rightarrow \eta_{U}(\hat{0})\,$ such that  $\,\{\eta_{V}(T(x_{n}))\}_{\Bbb{N}}=\{\hat{T}(\eta_{U}(x_{n}))\}_{\Bbb{N}}\rightarrow \eta_{V}(y)\,$.  Then there exists a sequence   Existe entonces una  $\,\{x_{n}'\}_{\Bbb{N}}$ in  $X\,$ such that  $\,\eta_{U}(x_{n}')=\eta_{U}(x_{n})\,$ and in addition,  $\,\{x_{n}'\}_{\Bbb{N}}\stackrel{\rho_{j}}{\rightarrow} \hat{0}\,$.  Therefore,  $\,\{\eta_{U}(x_{n}')\}_{\Bbb{N}}\rightarrow \eta_{U}(\hat{0})\,$.  We now obtain   $\,\{\hat{T}(\eta_{U}(x_{n}))\}_{\Bbb{N}}=\{\hat{T}(\eta_{U}(x_{n}'))\}_{\Bbb{N}}=\{(\hat{T}\circ \eta_{U})(x_{n}')\}_{\Bbb{N}}\rightarrow \eta_{V}(\tilde{0})\,$, which implies $\,\eta_{V}(y)\in Z_{\hat{\sigma}_{i}}\,$.\, By  Proposition 2.3,\, $\,\hat{T}\,$ is continuous.  \/ $\Box$    

\begin{cor}
Let  $\,U<X\,$\, y\, $\,V<Y\,$\, be closed vector subspaces  such that    $\,T(U)<V\,$.  Let   $\,\hat{T}:(X/U,\hat{\tau}) \rightarrow (Y/V,\hat{\sigma})\,$ be linear such that $\,\hat{T}\circ \eta_{U}= \eta_{V}\circ T\,$. (This is possible by the hypothesis).  Then:

\begin{itemize}
 \item[(i)]  If  $\,\hat{T}\,$ is continuous, then, $\,G(T)<V\,$. 
 \item[(ii)]  If for each  $\,i\in I\,$\, there exists $\,j\in J\,$ such that  $\,G_{ji}(T)<V\,$, then  $\,\hat{T}\,$ is continuous.
 \end{itemize}
\noindent  Here, $\,\eta_{U}:(X,\tau)\rightarrow (X/U,\hat{\tau})\,$ and $\,\eta_{V}:(Y,\sigma)\rightarrow (Y/V,\hat{\sigma})\,$ are the continuous canonical epimorphisms.   \end{cor}

\begin{defn}  Suppose $\{S_{n}:(X,\rho_{j_{n-1}})\rightarrow (S,\rho_{j_{n}})\}_{\Bbb{N}}\,$ and $\,\{R_{n}:(Y,\sigma_{i_{n}})\rightarrow (Y,\sigma_{i_{n-1}})\}_{\Bbb{N}}\,$\, are continuous linear functions such that     $\,T\circ S_{n}=R_{n}\circ T$  $\,\forall\, n\in \Bbb{N}\,$ and $\,\tau_{\rho_{n-1}}\subset \tau_{\rho_{n}}$.  Then we say that $\{S_{n}\}$ and  $\{R_{n}\}$ are {\bf intertwining}.  
\end{defn}

\noindent  The following is the main result.  It generalizes Lemma 1.6 of \cite{SC}.  In the proof we  will employ a gliding hump argument, applied to the seminorms.  \\

\begin{teo}
Let $\,\{S_{n}:(X,\rho_{j_{n-1}})\rightarrow (S,\rho_{j_{n}})\}_{\Bbb{N}}\,$ and $\,\{R_{n}:(Y,\sigma_{i_{n}})\rightarrow (Y,\sigma_{i_{n-1}})\}_{\Bbb{N}}\,$\, be continuous linear functions that are intertwining.   Then there exists  $\,N\in \Bbb{N}\,$ such that $\,\forall\, n\geq N\,$: 

 $$ \overline{R_{1}\circ R_{2}\circ ... \circ R_{n}}^{\sigma_{i_{o}}}=\overline{R_{1}\circ R_{2}\circ ... \circ R_{N}}^{\sigma_{i_{o}}}.   $$
\end{teo}

\noindent {\it Proof} 
   {\bf Claim 1}:  $\,\tau_{\rho_{j_{n-1}}}\subset \tau_{\rho_{n}}\, \Rightarrow\, G_{j_{n+1}i_{m}}(T)< G_{j_{n}i_{m}}(T) \; (\forall\, n, m\in \mathbb{N})$.    \\
 
\noindent  Indeed, for the continuity of the linear functions  $\,S_{n}:(X;\rho_{j_{n-1}})\rightarrow (X,\rho_{j_{n}})\,$\, there are positive numbers $\,\{c_{n}\}_{\Bbb{N}}\,$ such that   $\,\rho_{j_{n}}(S_{n}(x))\leq c_{n}\rho_{j_{n-1}}(x)\,$,  $\,\forall\, x\in X\,$.  We have  $\,y\in  G_{j_{n+1} i_{m}}(T)\,$ and $\,\{x_{n}\}_{\Bbb{N}}\stackrel{\rho_{j_{n+1}}}{\rightarrow} \hat{0}\,$ such that  $\,\{T(x_{n})\}_{\Bbb{N}}\stackrel{\sigma_{i_{m}}}{\rightarrow} y\,$ which tells us  $\,\{x_{n}\}_{\Bbb{N}}\stackrel{\rho_{j_{n}}}{\rightarrow} \hat{0}\,$ and  $\,\{T(x_{n})\}_{\Bbb{N}}\stackrel{\sigma_{i_{m}}}{\rightarrow} y\,$; that is,  $\,y\in G_{j_{n}i_{m}}(T)\,$.   \\

\noindent   {\bf Claim 2}:     For each  $\,n\in \Bbb{N}\,$,\, $\,R_{n+1}(G_{j_{n+1}i_{n+1}}(T))< G_{j_{n}i_{n}}(T)\,$ and  

$$\,\overline{R_{1}\circ R_{2}\circ ... \circ R_{n+1}}^{\sigma_{i_{o}}}=\overline{R_{1}\circ R_{2}\circ ... \circ R_{n}}^{\sigma_{i_{o}}}.   $$

\noindent By    Proposition 4.5-3,  and  Claim 1,   $\,R_{n+1}(G_{j_{n+1}i_{n+1}}(T))< G_{j_{n+1}i_{n}}(T)< G_{j_{n}i_{n}}(T)\,$.  From this,  $\,(R_{1}\circ R_{2}\circ ... \circ R_{n})(R_{n+1}(G_{j_{n+1}i_{n+1}}(T))< (R_{1}\circ R_{2}\circ ... \circ R_{n})(G_{j_{n}i_{n}}(T))\,$, and then   

$$\,\overline{R_{1}\circ R_{2}\circ ... \circ R_{n}\circ R_{n+1}(G_{j_{n+1}i_{n+1}}(T))}^{\sigma_{i_{o}}}<\overline{R_{1}\circ R_{2}\circ ... \circ R_{n}(G_{j_{n}i_{n}}(T))}^{\sigma_{i_{o}}}.  $$

\noindent  Let  $\,W_{n}=\overline{R_{1}\circ R_{2}\circ ... \circ R_{n}(G_{j_{n}i_{n}}(T))}^{\sigma_{i_{o}}}<Y\,$. We have that $\,W_{n+1}< W_{n}\,$.  For each  $\,n\in \Bbb{N}\,$\,  we consider the canonical linear epimorphism    $\,\eta_{n}:(Y,\sigma_{i_{o}})\rightarrow (Y/W_{n},\hat{\sigma}_{i_{o}})\,$, which is continuous.    \

\noindent   {\bf Claim 3}:  $\,\eta_{n}\circ R_{1}\circ ... \circ R_{n}\circ T: (X,\rho_{j_{n}})\rightarrow (Y/W_{n},\overline{\sigma}_{i_{o}})\,$\, is continuous.

\noindent     $\,((\eta_{n}\circ R_{1}\circ ... \circ R_{n})\circ T)(G_{j_{n}i_{n}}(T)=[\tilde{0}]_{W_{n}}=Z_{\overline{\sigma}_{i_{o}}}\,$.\,   By  Proposition 4.1,\, $\,\eta_{n}\circ R_{1}\circ ... \circ R_{n}\circ T\,$\, is continuous. From this, there exist positive numbers $\,\{d_{n}\}_{\Bbb{N}}\,$ such that   $\,\hat{\sigma}_{i_{o}}^{n}(\eta_{n}\circ R_{1}\circ ... \circ R_{n}\circ T(x))\leq d_{n} \rho_{j_{n}}(x)\,$ \hspace{2mm} $\,\forall\, x\in X\,$.    \

\noindent   {\bf Claim 4}:     $\,\eta_{n}\circ R_{1}\circ ... \circ R_{n-1}\circ T: (X,\rho_{j_{n-1}})\rightarrow (Y/W_{n},\overline{\sigma}_{i_{o}})\,$\, is not continuous if  $\,W_{n+1}\subsetneq W_{n}\,$ for infinitely many $n$.   \\
\noindent   If  $\,W_{n+1 }\subsetneq W_{n}\,$\hspace{2mm} for infinitely many  $ n\in \mathbb{N}\,$,\,   then by the discontinuity of the functions of Claim 4, we can construct the following sequence:  \
\noindent  Let   $\,x_{1}\in X\,$ such that  $\,\rho_{j_{1}}(x_{1})<\frac{1}{2}\,$ and $\,\hat{\sigma}_{i_{o}}^{2}((\eta_{2}\circ R_{1}\circ T)(x_{1}))> 1 + d_{2}\,$. Let  $\,x_{2}\in X\,$ such that $\,\rho_{j_{2}}(x_{2})<\frac{1}{2^{2}}\,$ and $\,\hat{\sigma}_{i_{o}}((\eta_{3}\circ R_{1}\circ R_{2}\circ T)(x_{2}))> 2 + d_{3} + \hat{\sigma}_{i_{o}}^{3}((\eta_{3}\circ T)(S_{1}(x_{1}))\,$.   Let     $\,x_{3}\in X\,$ such that  $\,\rho_{j_{3}}(x_{3})<\frac{1}{2^{3}}\,$ and $\,\hat{\sigma}_{i_{o}}((\eta_{4}\circ R_{1}\circ R_{2}\circ R_{3}\circ T)(x_{3}))> 3 + d_{4} + \hat{\sigma}_{i_{o}}^{4}(\eta_{4}\circ T)(S_{1}(x_{1})+(S_{2}\circ S_{1})(x_{2}))\,$.  In general, let   $\,x_{n}\in X\,$ such that $\,\rho_{j_{n}}(x_{n})<\frac{1}{2^{n}}\,$ and  

\begin{eqnarray} 
\hat{\sigma}_{i_{o}}((\eta_{n+1}\circ R_{1}\circ ...\circ R_{n}\circ T)(x_{n})) &> &n+ d_{n+1}  \nonumber \\ 
                                                                                                                    & + &\hat{\sigma}_{i_{o}}^{n+1}((\eta_{n+1}\circ T)(S_{1}(x_{1})+(S_{2}\circ S_{1})(x_{2}) \nonumber  \\
                                                                                                                    &+ &...+ (S_{n-1}\circ ... S_{1})(x_{n-1})))). \nonumber
\end{eqnarray} 

\noindent  Notice that for $\,m<n\in \Bbb{N}\,$,\, $\,\rho_{j_{m}} \leq c_{m+1}...c_{n}\rho_{j_{n}}\,$  and  \hspace{2mm} $\,\hat{\sigma}_{i_{o}}^{m}\leq \hat{\sigma}_{i_{o}}^{n}\,$.  For each  $\,n\in \Bbb{N}\,$ let  $\,z_{n}=(S_{1}\circ ... \circ S_{n})(x_{n})\in X\,$.   

\begin{eqnarray} 
\rho_{j_{o}}(z_{n})\leq c_{1}\rho_{j_{1}}((S_{2}\circ ... \circ S_{n})(x_{n})) &\leq & c_{1}c_{2}\rho_{j_{2}}((S_{3}   \circ  ... \circ S_{n})(x_{n}) \nonumber  \\
                                                                                                                   &  \leq & c_{1}c_{2} ... c_{n}\rho_{j_{n}}(x_{n})  \nonumber \\
                                                                                                                   & < &  (c_{1}c_{2}...c_{n})\frac{1}{2^{n}}   \nonumber   \end{eqnarray} 
\noindent  We consider two cases.\\
\noindent  {\bf Case 1}.  $\,c_{n}>1\,$\hspace{2mm} for infinitely many  $ n\in \mathbb{N}\,$.     \

For each $\,n\in \Bbb{N}\,$\, let the seminorm $\,\rho_{j_{n}}'=\frac{1}{c_{1}...c_{n}}\rho_{j_{n}}\,$.  Then  $\,\rho_{j_{n}}'(S_{n}(x))=\frac{1}{c_{1}...c_{n}}\rho_{j_{n}}(S_{n}(x))\leq \frac{1}{c_{1}...c_{n}}c_{n}\rho_{j_{n}}(x)=\rho_{j_{n-1}}'(x)\,$,  $\,\forall\, x\in X\,$.    $\,\overline{A}^{\rho_{j_{n}}}=\overline{A}^{\rho_{j_{n}}'}\,$ \hspace{2mm} $\,\forall\, A\subset X\,$.  The seminorms $\,\rho_{j_{n}}\,$\, and  $\,\rho_{j_{n}}'\,$ are equivalent.   \\

\noindent  {\bf Case 2}.   $\,c_{n}\leq 1\,$\ hspace{2mm} for infinitely many  $ n\in \mathbb{N}\,$.  

\noindent  We will first need to prove the following:  \\

\noindent   {\bf Claim 5}:    $\, \hat{\sigma}_{i_{o}}^{n+1}(\eta_{n+1}(T(\Sigma_{k=n+1}^{\infty}z_{k})))\leq \frac{1}{2^{n+1}}\,$.  \\

 To see this, we have  \begin{eqnarray} 
  \hat{\sigma}_{i_{o}}^{n+1}(\eta_{n+1}(T(\Sigma_{k=n+1}^{\infty}z_{k}))) & = &   \hat{\sigma}_{i_{o}}^{n+1}(\eta_{n+1}(T(S_{1}\circ ... \circ S_{n+1})( \Sigma_{k=2}^{\infty}(S_{n+2}\circ ... \circ S_{n+k})(x_{n+k}))))  \nonumber  \\
                                                                                                                   & = & \hat{\sigma}_{i_{o}}^{n+1}(\eta_{n+1}(R_{1}\circ ... \circ R_{n+1}\circ T)(\Sigma_{k=2}^{\infty}z_{k})))  \nonumber  \\
                                                                                                                   &  \leq &  d_{n+1} \rho_{j_{n}}(\Sigma_{k=2}^{\infty}(S_{n+2}\circ ... \circ S_{n+k})(x_{n+k}))   \nonumber  \\
                                                                                                                   &  \leq &  d_{n+1} \Sigma_{k=2}^{\infty}c_{n+2}...c_{n+k}\rho_{j_{o}}(x_{n+k})\leq d_{n+1}\Sigma_{k=2}^{\infty}\frac{1}{2^{n+k}}  \nonumber  \\
                                                                                                                   &= & d_{n+1}\frac{1}{2^{n+1}}.  \nonumber \end{eqnarray}

\noindent  Returning to the examination of Case  2, we have that   $\,\Sigma_{n=1}^{\infty}\rho_{i_{o}}(z_{n})\leq 1\,$.  As we required   $\,(X,\rho_{i_{o}})\,$ to be complete,  there exists   $\,z=\Sigma_{n=1}^{\infty}z_{n}\in X\,$.

\noindent     Otherwise,  $\,T(z_{n})=T((S_{1}\circ ... \circ S_{n})(x_{n}))=(R_{1}\circ ... \circ R_{n})(T(x_{n}))\,$.  \\

\noindent  Here comes the gliding hump part.  For each $\,n\in \Bbb{N}\,$\, we have:
\begin{eqnarray} 
\sigma_{i_{o}}(T(z))\geq \hat{\sigma}_{i_{o}}^{n+1}(\eta_{n+1}(T(z))) &  \geq  & \hat{\sigma}_{i_{o}}^{n+1}(\eta_{n+1}(T(\Sigma_{k=1}^{n}z_{k}))+\eta_{n+1}(T(\Sigma_{k=n+1}^{\infty}z_{k})) \nonumber \\
                                                                                                              & \geq &  \hat{\sigma}_{i_{o}}^{n+1}(\Sigma_{k=1}^{n}\eta_{n+1}((R_{1}\circ ... \circ R_{k})(T(x_{k}))) -  \nonumber  \\
                                                                                                              &  & \hat{\sigma}_{i_{o}}^{n+1}(\eta_{n+1}(T(\Sigma_{k=n+1}^{\infty}z_{k}))) \nonumber \\
                                                                                                              & \geq &\hat{\sigma}_{i_{o}}^{n+1}(\eta_{n+1}\circ R_{1} \circ ... \circ R_{n}\circ T(x_{n})) -  \nonumber  \\  
                                                                                                              & & \hat{\sigma}_{i_{o}}^{n+1}(\Sigma_{k=1}^{n-1}\eta_{n+1}((R_{1}\circ ... \circ R_{k})(T(x_{k}))) -  \nonumber  \\                                     
                                                                                                              &  &  \hat{\sigma}_{i_{o}}^{n+1}(\eta_{n+1}(T(\Sigma_{k=n+1}^{\infty}z_{k}))) \nonumber  \\
                                                                                                              & > &  n + d_{n+1} - \hat{\sigma}_{i_{o}}^{n+1}(\eta_{n+1}(T(\Sigma_{k=n+1}^{\infty}z_{k})))\geq n + d_{n+1} -  \nonumber  \\  
                                                                                                              &  & \frac{1}{2^{n+1}}d_{n+1} \nonumber  \\
                                                                                                              &  > & n  \nonumber . 
\end{eqnarray} 

\noindent   This contradiction  indicates to us that there can be only finitely many proper containments of the form   $\,W_{n+1}<W_{n}\,$.   \/ $\Box$   \\


\noindent  The following result generalizes  \cite{SC}; Lemma 1.8, page 13 to general locally convex spaces.   In the proof, we will have the pleasure of applying Zorn's Lemma several times.  \\
\begin{teo}
Let     $\,\{R_{n}:(Y,\sigma)\rightarrow (Y,\sigma)\}_{n\in \Bbb{N}}\,$\, be a sequence of commuting continuous linear functions.  Then the following properties hold: \

\begin{enumerate}
\item   There exists  $\,Y_{\infty}<Y\,$ such that  $\,R_{n}(Y_{\infty})=Y_{\infty}\,$\, $\,\forall\, n\in \Bbb{N}\,$. 
\item  If  $\,V<Y\,$ is such that  $\,R_{n}(V)=V\,$\, $\,\forall\, n\in \Bbb{N}\,$,\, then \, $\,V<Y_{\infty}\,$. 
\item  There exists  $\,Y^{\infty}<Y\,$\, $\,\sigma$ - closed such that  $\,\overline{R_{n}(Y^{\infty})}^{\sigma}=Y^{\infty}\,$\,  $\,\forall\, n\in \Bbb{N}\,$. 
\item  If  $\,W<Y\,$ is  $\,\sigma$ -  closed and satisfies    $\,\overline{R_{n}(W)}^{\sigma}=W\,$\,  $\,\forall\, n\in \Bbb{N}\,$\, then  $\,W<Y^{\infty}\,$. 
\item  $\, \overline{R_{n}(Y^{\infty})}^{\sigma}=\overline{R_{n}(R_{m}(Y^{\infty}))}^{\sigma}=Y^{\infty}\,$ \hspace{2mm} $\,\forall\,n,\,m\in \Bbb{N}\,$.  
\item  $\,\overline{Y_{\infty}}^{\sigma}<Y^{\infty}\,$. \hspace{6cm}      
\end{enumerate}  
\end{teo}

\noindent {\it Proof} 
  1.   Let us  say that $\,V<Y\,$\, satisfies {\it property} $P (\models P(V))$, if   $\,R_{n}(V)=V\,$\,  $\,\forall\, n\in \Bbb{N}\,$.  \
Let   $\,F=\{V<Y\,|\,\models P(V)\}\,$\,  be ordered by inclusion. Observe that $\,\{\tilde{0}\}\in F\,$. Also, let $\,C=\{V_{\lambda}\}_{\Lambda}\subset F\,$\, be a chain.  Then $\,V=\bigcup_{\Lambda}V_{\lambda}<Y\,$.   On the other hand, for each  $\,n\in \Bbb{N}\,$,\, $\,R_{n}(V)=\bigcup_{\Lambda}R_{n}(V_{\lambda})=\bigcup_{\Lambda}V_{\lambda}=V\,$.  In other words,  $\,V\in F\,$\,  and is a maximal element of  $\,C\,$.  Zorn's Lemma now guarantees the existence of a maximal subspace in $\,F\,$.   Denote said maximal subspace by  $\,Y_{\infty} \in F$.  Now, in item 2, we will prove that $\,Y_{\infty}$ is unique.   \\

\noindent  2.  We will prove that if  $\,V\in F\,$, then \, $\,V+Y_{\infty}\in F\,$: To this end, we have $\,R_{n}(V+Y_{\infty})=R_{n}(V)+R_{n}(Y_{\infty})=V+Y_{\infty}\,$.   $\,Y_{\infty}<V+Y_{\infty}\, \Rightarrow\, Y_{\infty}=V+Y_{\infty}\, \Rightarrow\, V<Y_{\infty}\,$. 

\noindent  3.  For this part, let us say that  $\,W<Y\,$\, satisfies  {\it property}\, $\,Q\,$\, ($\,\models Q(W)\,$)\, if  $\,\overline{R_{n}(W)}^{\sigma}=W\,$\,  $\,\forall\, n\in \Bbb{N}\,$.  Now, let  $\,G=\{W<Y\,|\,\models Q(W)\}\,$\, be ordered by inclusion.   and observe that  $\,\{\tilde{0}\}\in G\,$.   Let  $\,D=\{W_{\lambda}\}_{\Lambda}\subset G\,$\,be a chain.  We have that  $\,W=\overline{\bigcup_{\Lambda}W_{\lambda}}<Y\,$. For each  $\,n\in \Bbb{N}\,$: 

 \begin{eqnarray} 
 R_{n}(W)<\overline{R_{n}(\bigcup_{\Lambda}W_{\lambda})}^{\sigma} &=& \overline{\bigcup_{\Lambda}R_{n}(W_{\lambda})}^{\sigma}    \nonumber  \\
                                                                                                                     &<& \overline{\bigcup_{\Lambda}\overline{R_{n}(W_{\lambda})}}^{\sigma}    \nonumber  \\
                                                                                                                     &=&  \overline{\bigcup_{\Lambda}W_{\lambda}}^{\sigma}   \nonumber  \\
                                                                                                                     &=&W.     \nonumber  \end{eqnarray}   

\noindent  From the above we obtain $\,\overline{R_{n}(W)}^{\sigma}<W\,$.  Let   $\,y\in W\,$\,  and let  $\,O\subset Y\,$\, be a closed, absolutely convex neighborhood of  $\,\tilde{0}\,$.  Finally, let  $\,O'\,$\,be another neighborhood of  $\,\tilde{0}\,$\, for which  $\,O'+O'\subset O\,$.  By the  $\,\sigma$ - density of the subspace $\,\bigcup_{\Lambda}W_{\lambda}\,$ in $\,W\,$, there exists  $\,w\in W_{\lambda}\,$\, (for some $\,\lambda\in \Lambda\,$) such that  $\,w\in y+O'\,$.   Meanwhile, by the $\,\sigma$- density of $\,R_{n}(W_{\lambda})\,$\, in\, $\,W_{\lambda}\,$,\, there exists $\,w'\in W_{\lambda}\,$\,  such that   $\,R_{n}(w')\in w+O'\,$.  We conclude from this that $\,R_{n}(w')\in w+O'\subset (y+O')+O'\subset y+O\,$.  In other words,  $\,y\in \overline{R_{n}(W)}^{\sigma}\,$.  We have shown that  $\,W\in G\,$\, and is an upper bound for  $\,D\,$.  Here comes Zorn's Lemma again:  $\,G\,$\, has a maximal element.  Put  $\,Y^{\infty}\in G\,$ as a maximal element, and in the next item we will prove that it is unique.    

\vspace{2mm}

\noindent  4.  We will prove that if    $\,W\in G\,$ then $\,\overline{W+Y^{\infty}}^{\sigma}\in G\,$: \
For one direction we have:  \begin{eqnarray} 
R_{n}(\overline{W+Y^{\infty}})& <&  \overline{R_{n}(W)+R_{n}(Y^{\infty})}^{\sigma}    \nonumber  \\
                                                      & < & \overline{R_{n}(W)}^{\sigma}+ \overline{R_{n}(Y^{\infty})}^{\sigma}    \nonumber  \\
                                                      &  < & \overline{\overline{R_{n}(W)}^{\sigma}+ \overline{Y^{\infty}}^{\sigma}}^{\sigma}    \nonumber  \\
                                                      & =&  \overline{W+Y^{\infty}}^{\sigma}.   \nonumber \end{eqnarray}      \\

\noindent   Now we will prove the other direction:  Let   $\,y\in \overline{W+Y^{\infty}}^{\sigma}\,$ and let  $\,O\subset Y\,$\, be a closed, absolutely convex zero neighborhood.  Also, let  $V$ be another zero neighborhood for which $\,V+V\subset O\,$.  There exist, then,  $\,w\in W=\overline{R_{n}(W)} ^{\sigma}\,$ and  $\,z\in Y^{\infty}=\overline{R_{n}(Y^{\infty})}^{\sigma}\,$ such that $(w+z)\in y+V\,$.  Moreover, there exist  $\,w'\in W\,$ and  $\,z'\in Y^{\infty}$ such that  $\,R_{n}(w')\in w+V$ and  $\,R_{n}(z')\in z+V$.  We  have  
 $$\,R_{n}(w'+z')=R_{n}(w')+R_{n}(z')\in (w+ V)+(z+V) = (w+z)+(V+V)\subset y+O.  $$

\noindent  In other words,   $\,y\in \overline{R_{n}(W+Y^{\infty})}^{\infty}\,$.  Thus, we have proved that  $\,\overline{W+Y^{\infty}}^{\sigma}\in G\,$.  Now, 

$$\,Y^{\infty}=\overline{Y^{\infty}}^{\sigma}<\overline{W+Y^{\infty}}^{\sigma}\, \Rightarrow\, Y^{\infty}= \overline{W+Y^{\sigma}}^{\sigma}\, \Rightarrow\, W<Y^{\infty}.  $$

\noindent  5. Observe:
  $$\,\overline{R_{n}(R_{m}(Y^{\infty}))}^{\sigma}> R_{n}(\overline{R_{m}(Y^{\infty})}^{\sigma})=R_{n}(Y^{\infty}).  $$

\noindent    This implies that  $$\,Y^{\infty}=\overline{R_{n}(Y^{\infty})}^{\sigma}<\overline{R_{n}(R_{m}(Y^{\infty}))}^{\sigma}<\overline{R_{n}(Y^{\infty})}^{\sigma}=Y^{\infty}.  $$

\noindent     On the other hand,    $$\,Y^{\infty}=\overline{R_{n}(Y^{\infty})}^{\sigma}=\overline{R_{n}(R_{m}(Y^{\infty}))}^{\sigma}=\overline{R_{n}(Y^{\infty})}^{\sigma}.  $$

\vspace{2mm}

\noindent  6. First, let us recall that for each  $\,A\subset Y\,$,\, $\,\overline{A}^{\sigma}=\bigcap_{I}\overline{A}^{\sigma_{i}}\,$.  \\

\noindent   {\bf Claim}:    $\,\overline{Y_{\infty}}^{\sigma}\in G\,$.  To see this, let 

$$\,y\in \overline{R_{n}(\overline{Y_{\infty}}^{\sigma})}^{\sigma}=\bigcap_{I}\overline{R_{n}(\overline{Y_{\infty}})}^{\sigma_{i}}.  $$

  Thus, for each     $\,i\in I\,$ as well as each $\,\epsilon >0\,$, we have:  There exists $\,z_{i}\in \overline{Y_{\infty}}^{\sigma}\,$ such that  $\,\sigma_{i}(y-z_{i})<\frac{\epsilon}{2}\,$\,  and there is  $\,y_{i}\in Y_{\infty}\,$ such that $\,\sigma_{i}(y_{i}-z_{i})<\frac{\epsilon}{2}\,$.  Hence, $\,\sigma_{i}(y-y_{i})<\epsilon\,$.  This proves that $\,y\in \overline{Y_{\infty}}^{\sigma_{i}}\,$ for each  $\,i\in I\,$\, and therefore  $\,y\in \overline{Y_{\infty}}^{\sigma}\,$.This last statement tells us that  $\,\overline{R_{n}(\overline{Y_{\infty}}^{\sigma})}<\overline{Y_{\infty}}^{\sigma}\,$.  Now we prove the other direction:   $R_{n}(Y_{\infty})  = Y_{\infty}\,\Rightarrow $ 

$$ \overline{Y_{\infty}}^{\sigma}    = \overline{R_{n}(Y_{\infty})}^{\sigma}   = \bigcap_{I}\overline{R_{n}(Y_{\infty})}^{\sigma_{i}}   <   \bigcap_{I}\overline{R_{n}(\overline{Y_{\infty}}^{\sigma})}^{\sigma_{i}}=\overline{R_{n}(\overline{Y_{\infty}}^{\sigma})}^{\sigma}.  $$

\noindent  Putting both directions together proves that  $\,\overline{Y_{\infty}}^{\sigma}\in G\,$ and as a consequence,  $\,\overline{Y_{\infty}}^{\sigma}<Y^{\infty}. \; \; \Box$ \\

\noindent {\bf Acknowledgements.}  The authors C. Bosch, C. Garc\'{i}a, C. G\'{o}mez-Wulschner, and R. Vera   were partially supported by the Asociaci\'{o}n Mexicana de Cultura, A.C.   \\


\bigskip

\begin{thebibliography}{99}

\bibitem{HQ}    He, F.,   J. Qiu, J.:  Sequentially lower complete spaces and Ekeland's variational principle, {\it Acta Math. Sin. (Engl. Ser.)},  {\bf 31},  no. 8, 1289 - 1302 (2015).     \\

\bibitem{NB}      Narici, L.,   Beckenstein, E.:  Topological Vector Spaces, M. Dekker, New York,  (1985).  \\

\bibitem{NA}    Neumann, M.,  Albrecht, E.:   Automatic continuity of generalized local linear operators, {\it Manuscripta Math.}, {\bf 32},  263 - 294 (1980).    \\

\bibitem{PC}    P. P\'{e}rez Carreras, P.,   Bonet, J.:  Barrelled Locally Convex Spaces, North - Holland \# {131}, Amsterdam,  (1987).  \\


\bibitem{QH}    Qiu, J.,   He, F.:  p-distances, q-distances and a generalized Ekeland's variational principle in uniform spaces, {\it  Acta Math. Sin. (Engl. Ser.)},  {\bf 28}, no. 2, 235 - 254  (2012).   \\


\bibitem{SC}   Sinclair, A.:  Automatic Continuity of Linear Operators, Cambridge University Press, Cambridge, (1976).  Online edition published by Cambridge University Press, 20 January, 2009.  \\

\bibitem{CS}  Swartz, C.:  Infinite Matrices and the Gliding Hump, World Scientific, Singapore, (1996).  \\

\bibitem{AT}     Tajmouati, A.: On automatic boundedness of linear operators on convex bornological spaces, {\it Ital. J. Pure \& Applied Math.}, {\bf 32}, 155 - 164  (2014).   

\end{thebibliography}
\end{document}